\documentclass[11pt]{amsart}

\usepackage[T1]{fontenc}
\usepackage{amsmath,amssymb,amsthm,mathtools}
\usepackage{enumitem}
\usepackage{booktabs}
\usepackage{array}
\usepackage{longtable}
\usepackage{cite}
\usepackage[colorlinks=true,linkcolor=blue,citecolor=blue,urlcolor=blue]{hyperref}
\hypersetup{bookmarksnumbered=true,bookmarksopen=true,bookmarksopenlevel=2}
\newtheorem{theorem}{Theorem}[section]
\newtheorem{proposition}[theorem]{Proposition}
\newtheorem{lemma}[theorem]{Lemma}
\newtheorem{corollary}[theorem]{Corollary}
\theoremstyle{definition}
\newtheorem{definition}[theorem]{Definition}
\theoremstyle{remark}
\newtheorem{remark}[theorem]{Remark}
\newcommand{\K}{\mathbb K}
\newcommand{\R}{\mathbb R}
\newcommand{\C}{\mathbb C}
\newcommand{\N}{\mathbb N}
\newcommand{\F}{\mathcal F}
\newcommand{\B}{\mathcal B}
\newcommand{\Prob}{\operatorname{Prob}}
\newcommand{\supp}{\operatorname{supp}}
\newcommand{\ran}{\operatorname{ran}}
\newcommand{\spn}{\operatorname{span}}
\newcommand{\Id}{\operatorname{Id}}
\newcommand{\Rea}{\operatorname{Re}}
\newcommand{\norm}[1]{\left\lVert #1\right\rVert}
\newcommand{\abs}[1]{\left\lvert #1\right\rvert}
\newcommand{\tens}{\widehat{\otimes}_{\pi}}
\numberwithin{equation}{section}

\title[Approximation property in dual Banach spaces]{Dual Banach spaces with AP but without BAP}

	\author[Longfa Sun]{Longfa Sun}
\address[L. Sun]{Hebei key laboratory of physics and energy technology, School of mathematics and physics, North China Electric Power University, Baoding 071003, China}
\email{sun.longfa@ncepu.edu.cn}

\author[Zhenguo Wei]{Zhenguo Wei}
\address[Z. Wei]{School of mathematics and statistics, Wuhan university, 430070, Wuhan, Hubei, China}
\email{zhenguo.wei@whu.edu.cn}

\author[Zhendong  Xu]{Zhendong Xu}
\address[Z. Xu]{%Institute for Advanced Study in Mathematics, Harbin Institute of Technology,  Harbin 150001, China; and
	Department of Mathematical Sciences and the Research Institute of Mathematics, Seoul National University, Gwanak-ro 1, Gwanak-gu, Seoul 08826, Republic of Korea}
\email{zhendong\_xu\_97@snu.ac.kr}

\author[Hao Zhang]{Hao Zhang}
\address[H. Zhang]{%Institute for Advanced Study in Mathematics, Harbin Institute of Technology,  Harbin 150001, China; and
	Instituto de Ciencias Matem\'{a}ticas, Consejo Superior de Investigaciones Cient\'{i}ficas, C/ Nicol\'{a}s Cabrera 13-15, 28049, Madrid, Spain}
\email{hao.zhang@icmat.es}

\thanks{
Longfa Sun is supported by the Fundamental Research Funds
for the Central Universities (Grant No. 2025MS178). Zhendong Xu is partially supported by the National Research Foundation of Korea (NRF) grant funded by the Korea government (MSIT) (No.2020R1C1C1A01009681), Samsung Science and Technology Foundation under Project Number SSTF-BA2002-01 and the Basic Science Research Program through the National Research Foundation of Korea (NRF) (Grant RS-2022-NR069971).
}

\date{}

\date{}


\begin{document}
\begin{abstract}
We give negative answers to the AP-to-BAP and BAP-to-MAP questions for dual Banach spaces. More precisely, there is a separable Banach space $V$ with the metric approximation property (MAP) such that $V^*$ has the approximation property (AP) but no bounded approximation property (BAP).
There is also a separable MAP space $X$, isomorphic to $\ell_1$, whose dual has BAP but no MAP. This resolves a longstanding open question explicitly stated in \cite{Oja2010}. The constructions rely on a finitely generated group without Property A and a quantitative estimate for optimal bounded approximation constants of projective tensor products.
\end{abstract}

\maketitle
\tableofcontents

\section{Introduction}

The approximation property (AP) is one of the fundamental concepts in the geometry of Banach spaces and the theory of linear operators. Its origins can be traced back to the early development of functional analysis, particularly to Banach's monograph \cite{Banach1932} and the approximation problems raised by Mazur in the Scottish Book \cite{ScottishBook}. The first systematic investigation of the approximation property was undertaken by Grothendieck in his celebrated memoir \cite{Grothendieck}. His work established a profound connection between finite rank approximation, tensor products, nuclear operators and the structure of Banach spaces. In particular, Grothendieck obtained several equivalent characterizations of the approximation property and introduced its metric version, thereby laying the foundations for much of the subsequent development of approximation theory.

A Banach space $X$ has the \emph{approximation property} (AP), if for every norm-compact set $K\subset X$ and every $\varepsilon>0$, there exists a bounded finite rank operator $T:X\to X$ such that
\[
\sup_{x\in K}\norm{Tx-x}<\varepsilon.
\]
For $\lambda\geq1$, the $\lambda$-\emph{bounded approximation property} ($\lambda$-BAP) requires $\norm{T}\leq\lambda$. The BAP means $\lambda$-BAP for some finite $\lambda$, and the \emph{metric approximation property} (MAP) is $1$-BAP.
Consequently,
\[
\mathrm{MAP}\Longrightarrow\mathrm{BAP}\Longrightarrow\mathrm{AP}.
\]
For a nonzero Banach space $X$, define its optimal bounded approximation constant by
\begin{equation}\label{eq:beta}
	\beta(X)=\inf\{\lambda\geq1:X\text{ has }\lambda\text{-BAP}\},
	\quad \inf\varnothing=\infty.
\end{equation}

It was unknown whether every Banach space possesses the approximation property for a long time. This fundamental problem was resolved negatively by Enflo \cite{Enflo1973} in 1973, through the construction of a separable reflexive Banach space without AP. In the same year, Figiel and Johnson \cite{FigielJohnson} demonstrated that AP does not imply BAP, thereby establishing a genuine distinction between qualitative finite rank approximation and approximation under a uniform norm constraint. These discoveries stimulated extensive research on approximation properties, their quantitative constants and their relationships with the local structure of Banach spaces. We refer the interested reader to Casazza's survey \cite{Casazza} for a comprehensive account of the classical theory.

For the theory of operator algebras, one of the most important developments of approximation properties is the characterization of nuclear $C^*$-algebras in terms of finite-dimensional completely positive approximations; see Choi and Effros \cite{ChoiEffros1978}. This provides a fundamental connection between approximation properties and tensor product structures of operator algebras. However, Haagerup \cite{Haagerup1979} found that the reduced $C^*$-algebras of free groups are non-nuclear but have MAP, which distinguishes MAP from nuclearity even within the class of $C^*$-algebras. On the negative side, Szankowski \cite{Szankowski1981} showed that $B(\ell_2)$ fails AP.
\iffalse
Approximation properties also play a central role in the theory of operator algebras. One of the most important developments is the characterization of nuclear $C^*$-algebras in terms of finite-dimensional completely positive approximations, established through the work of Choi and Effros \cite{ChoiEffros1978}. This provides a fundamental connection between approximation properties and the tensor product structure of operator algebras. In another direction, Haagerup \cite{Haagerup1979} proved that the reduced $C^*$-algebra of a free group has the metric approximation property, despite being nonnuclear. Thus, even within the class of $C^*$-algebras, metric approximation and nuclearity are distinct notions. On the negative side, Szankowski \cite{Szankowski1981} showed that the algebra of all bounded operators on an infinite-dimensional separable Hilbert space fails the approximation property. These results illustrate the subtle behavior of finite-rank approximation in operator algebras.
\fi
The development of operator space theory led to further refinements of classical approximation properties. Replacing bounded linear maps by completely bounded maps and taking into account matrix norm structures of operator spaces, one can define the completely bounded approximation property (CBAP) and the operator space approximation property (OAP). These notions were systematically investigated by Effros and Ruan \cite{EffrosRuan1990}; see also \cite{EffrosRuan2000}.
For group operator algebras, these properties are closely related to weak amenability~\cite{CowlingHaagerup1989} and the Haagerup--Kraus
approximation property~\cite{HaagerupKraus1994}. These developments revealed a close relationship between harmonic analysis on groups and the geometry of operator spaces.

\iffalse
The development of operator space theory led to further refinements of classical approximation properties. In this framework, completely bounded maps replace bounded linear maps, giving rise to the completely bounded approximation property (CBAP) and the operator space approximation property (OAP). These notions were systematically investigated by Effros and Ruan \cite{EffrosRuan1990}; see also \cite{EffrosRuan2000}. They provide a natural framework for studying finite-dimensional approximation while taking into account the matrix norm structure of operator spaces. In the setting of group operator algebras, Cowling and Haagerup \cite{CowlingHaagerup1989} developed the theory of weak amenability and completely bounded Fourier multipliers, while Haagerup and Kraus \cite{HaagerupKraus1994} introduced an approximation property for locally compact groups and established important connections with approximation properties of their associated $C^*$-algebras and von Neumann algebras. These developments revealed a close relationship between harmonic analysis on groups and the geometry of operator spaces.
\fi

Approximation properties have also become an important tool in the study of noncommutative $L_p$-spaces. Junge and Ruan \cite{JungeRuan2003}
% investigated approximation properties of noncommutative $L_p$-spaces associated with discrete groups. They
showed that the Haagerup--Kraus approximation property of a discrete group implies the OAP of the corresponding noncommutative $L_p$-spaces for $1<p<\infty$. Under additional hypotheses, they obtained stronger approximation properties and consequences for the existence of completely bounded Schauder bases. Further developments by Junge, Nielsen, Ruan and Xu \cite{JungeNielsenRuanXu2004} connected approximation properties with the local structures of noncommutative $L_p$-spaces and finite-dimensional Schatten classes. In contrast, Lafforgue and de la Salle \cite{LafforgueDeLaSalle2011} constructed noncommutative $L_p$-spaces without the CBAP, using harmonic analysis and completely bounded Schur multipliers on higher-rank groups. Related results were subsequently obtained by de Laat \cite{DeLaat2013}. These investigations show that approximation properties are not merely technical tools, but can reflect deep geometric and analytic features of the underlying operator algebras.

Despite these extensive developments, several basic questions concerning the classical approximation properties of dual Banach spaces have remained unresolved. The distinction between AP, BAP, and MAP is particularly subtle in this setting. Grothendieck's theory \cite{Grothendieck} implies that a separable dual Banach space with AP necessarily has MAP \cite[Theorem 1.e.15]{LT}. More generally, positive results are available under additional geometric assumptions, including the Radon--Nikodym property (RNP). More precisely, if $X^*$ has AP and either $X^*$ or $X^{**}$ has the RNP, then $X^*$ has MAP \cite[Theorem 3.2]{Oja2010}. However, these results do not apply to arbitrary dual spaces. The question of whether AP of a dual Banach space always implies MAP, already rooted in Grothendieck's work, has been discussed extensively in the literature; for instance, the question for general duals is listed in Casazza's survey \cite{Casazza} and explicitly stated by Oja \cite[Problem 2.1]{Oja2010}; see also the discussions in \cite{Oja2014}. Oja further separates it into two subquestions on Page 221. We formulate these questions as follows:
\begin{enumerate}[label=(Q\arabic*),ref=Q\arabic*,leftmargin=*]
	\item\label{q:ap-map} Does AP of $X^*$ imply MAP of $X^*$?
	\item\label{q:ap-bap} Does AP of $X^*$ imply BAP of $X^*$?
	\item\label{q:bap-map} Does BAP of $X^*$ imply MAP of $X^*$?
\end{enumerate}
All questions above are left open.

The main purpose of the present paper is to give negative answers to all questions. We would like to remark that AP and BAP are invariant under isomorphisms, while MAP depends on the specified norm. Throughout this paper, a dual space is equipped with the induced dual norm.

\medskip
% \subsection{Main results}
Our first main result answers \eqref{q:ap-bap}, and therefore \eqref{q:ap-map}, negatively.
\begin{theorem}\label{thm:main-ap}
	Over either $\R$ or $\C$, there exists a separable Banach space $X$ with MAP such that $X^*$ has AP and does not have BAP.
\end{theorem}

%Theorem~\ref{thm:main-ap} has another consequence for a question discussed in \cite[Section 3]{Oja2010}. By Lima and Oja \cite{LimaOja}, AP and weak MAP are equivalent for dual spaces. Our example therefore has weak MAP but no BAP. This separates weak MAP from BAP even for a dual with a separable MAP predual. The definitions and the precise corollary appear in Section~\ref{sec:tensor}.

Our second result answers \eqref{q:bap-map} negatively within the isomorphism class of $\ell_\infty$, using a norm that is explicitly induced by a separable predual.
\begin{theorem}\label{thm:main-metric}
	Over either $\R$ or $\C$, for every $0<a\leq 1$, there is a separable Banach space $X$, isomorphic to $\ell_1$ and having MAP, such that $X^*$ has BAP and
	\[
	1<\beta(X^*)\leq1+a.
	\]
	Consequently, $X^*$ does not have MAP.
	\iffalse
	There are a finite set $F\subset X^*$ and $\gamma>0$ such that every finite-rank contraction $H:X^*\to X^*$ satisfies
	\[
	\max_{f\in F}\norm{Hf-f}\geq\gamma.
	\]
	\fi
\end{theorem}
Thus the failure of MAP persists even when the optimal bounded approximation constant is arbitrarily close to $1$.

We would like to mention that our results also have direct consequences for weaker forms of approximation property. Lima and Oja \cite{LimaOja} introduced the weak metric approximation property (weak MAP) and showed that, for dual Banach spaces, it is equivalent to the AP. Consequently, our construction provides a dual Banach space with weak MAP but without BAP. (The definitions of weak MAP and the precise statement appear at the end of Section~\ref{sec:tensor}.) This further illustrates the distinction between qualitative approximation, bounded approximation and their metric variants, even within the class of dual Banach spaces.

The strategy of our constructions combines ideas from the geometry of Banach spaces, geometric group theory and the theory of tensor products. The starting point is a finitely generated group without Yu's Property A \cite{Yu}. The first examples of finitely generated groups without Property A were constructed by Gromov \cite{Gro03}; see also other examples of finitely generated groups without Property A constructed by Osajda \cite{Osajda,Osa18}. By encoding the group action into an equivalent norm on a space isomorphic to $\ell_1$, we obtain a dual Banach space whose optimal bounded approximation constant admits a quantitative lower bound in terms of an obstruction to Property A. Then we establish a multiplicative lower estimate for the optimal bounded approximation constants of duals of projective tensor products. By iterating this estimate, we obtain a sequence of dual spaces with finite but arbitrarily large approximation constants. Finally, by applying a suitable direct-sum construction introduced by \cite{FigielJohnson}, we construct a dual Banach space with AP but without BAP, while preserving MAP and separability of its predual.

%Can finite rank approximation of the identity on a dual Banach space always be achieved with a uniform bound on the approximating operators? If a bound is available, can it be reduced to $1$ without changing the dual norm? These questions distinguish three classical approximation properties. This paper constructs dual spaces giving negative answers to both questions, while their preduals are separable and admit contractive finite rank approximations.

% The finite set in the theorem witnesses an obstruction to \emph{every} finite-rank contraction.

% No positivity, weak-star continuity, or coordinate form is imposed on $H$. The equivalence $\beta(E)=1\Longleftrightarrow\mathrm{MAP}$ is justified in Lemma~\ref{lem:finite-tests} below.

This paper is organized as follows. Section~\ref{sec:preliminaries} is devoted to the preliminaries on the notions used in the proofs. In Sections~\ref{sec:construction}, we prove Theorem~\ref{thm:main-metric}. In the last section, we show Theorem~\ref{thm:main-ap} and the non-equivalence between weak MAP and BAP.

\

\noindent \textbf{Notation:} The scalar field is \(\K\in\{\R,\C\}\). Denote by \(\F(E,F)\) the set of all bounded finite rank operators from \(E\) to \(F\), and \(\F(E)=\F(E,E)\). For a linear operator $T$, denote its range by $\ran T$.
% The canonical embedding is \(J_E:E\to E^{**}\).
On \(\ell_1(I)\), we use its usual standard basis \((e_i)_{i\in I}\). For finite \(J\subset I\), \(P_J\) is the canonical coordinate projection. Equivalent norms will always be identified explicitly.

\bigskip

\section{Preliminaries}\label{sec:preliminaries}

This section is devoted to a brief introduction to finite tests and Property A. In this section, $X$ always denotes a Banach space. We proceed with finite tests.

\subsection{Finite tests}
We give a self-contained explanation of using finite tests to describe BAP; the interested reader is referred to \cite{Grothendieck,LT,Casazza} for more details. Let $X\neq\{0\}$ be a Banach space. A \emph{finite test} on $\B(X)$, denoted by $\ell$, is a functional such that there exist $m\in \mathbb{N}$, $f_j\in X,\ g_j\in X^*$ $(1\leq j\leq m)$ and
\begin{equation}\label{h:equa2.0}
 \ell(T)=\sum_{j=1}^m g_j(Tf_j),\quad  \forall  T \in \B(X).
\end{equation}
Let $\ell$ be a finite test, and define
\begin{equation}\label{eq:test-constants}
 c_X(\ell)=\sup_{T\in\F(X),\,\norm{T}\leq1}\abs{\ell(T)},
 \quad t_X(\ell)=\ell(\Id_X).
\end{equation}
The first quantity $c_X(\ell)$ is finite because it is at most $\sum_j\norm{g_j}\norm{f_j}$.

\begin{lemma}\label{lem:finite-tests}
Let $X\neq\{0\}$ be a Banach space and $\lambda\geq1$. Then $X$ has $\lambda$-BAP if and only if for every finite test $\ell$,
\begin{equation}\label{h:equa2.1}
 \abs{t_X(\ell)}\leq\lambda c_X(\ell).
\end{equation}
Moreover,
\begin{equation}\label{eq:test-beta}
 \beta(X)=\sup_{c_X(\ell)>0; \ \ell \ \text{finite test}}\frac{\abs{t_X(\ell)}}{c_X(\ell)}.
\end{equation}
If $\beta(X)<\infty$, then $X$ has $\beta(X)$-BAP. In particular, $\beta(X)=1$ if and only if $X$ has MAP.
\end{lemma}
\begin{proof}
We first show \eqref{h:equa2.1}. Suppose that $X$ has $\lambda$-BAP, and $\ell$ is a finite test on $\B(X)$. Then there exist $m\in \mathbb{N}$, $f_j\in X,\ g_j\in X^*$ $(1\leq j\leq m)$ and
\[
\ell(T)=\sum_{j=1}^m g_j(Tf_j),\quad  \forall  T \in \B(X).
\]
Given $f_1, \cdots, f_m \in X$, there exists a sequence of finite rank operators $\{T_n\}_{n \geq 1} \subset \F(X)$ such that
\[
\| T_n \| \leq \lambda, \quad \text{and} \ \max_{1 \leq i \leq m} \| T_nf_i - f_i \| \to 0, \quad \text{ as } n \to \infty.
\]
Then $\ell(T_n)\to t_X(\ell)$ as $n \to \infty$, and this implies
 $$ |t_X(\ell)|=\lim_{n\to \infty} |\ell(T_n)| = \lim_{n\to \infty}\lambda | \ell ( \lambda^{-1}T_n )| \leq \lambda c_X(\ell) .$$
This yields \eqref{h:equa2.1}.

Conversely, suppose that \eqref{h:equa2.1} holds. We need to show that $X$ has $\lambda$-BAP. Fix $m\in \mathbb{N}$ and $f_1,\ldots,f_m\in X$.  Let $X^m$ be equipped with its maximum norm. Then the following set
\[
 \mathcal C_\lambda=\{(Tf_1,\ldots,Tf_m):T\in\F(X),\ \norm{T}\leq\lambda\}
\]
is convex and balanced. Therefore, $X$ has $\lambda$-BAP if any fixed $(f_1,\ldots,f_m)\in X^m$ belongs to the norm closure of $\mathcal{C}_\lambda$. Indeed, it suffices to show that approximation on finite subsets implies uniform approximation on norm-compact subsets. Let $K\subset X$ be a nonempty norm-compact set and $\varepsilon>0$. Since $K$ is compact, there exist $f_1,\ldots,f_m\in K$ such that
\[
K\subseteq\bigcup_{j=1}^m
B\left(f_j,\frac{\varepsilon}{3(\lambda+1)}\right).
\]
By the preceding argument, we have $(f_1,\ldots,f_m)\in\overline{\mathcal C_\lambda}.$ Hence there exists $T\in\mathcal F(X)$ with $\|T\|\leq\lambda$ such that
\[
\max_{1\leq j\leq m}\|Tf_j-f_j\|
<\frac{\varepsilon}{3}.
\]
For any $x\in K$, choose $1\leq j\leq m$
such that
\[
\|x-f_j\|<\frac{\varepsilon}{3(\lambda+1)}.
\]
Then by the triangle inequality,
\[
\begin{aligned}
	\|Tx-x\|
	&\leq \|Tx-Tf_j\|
	+\|Tf_j-f_j\|+\|f_j-x\|\\
	&\leq (\|T\|+1)\|x-f_j\|
	+\|Tf_j-f_j\|\\
	&<(\lambda+1)
	\frac{\varepsilon}{3(\lambda+1)}
	+\frac{\varepsilon}{3}\\
	&=\frac{2\varepsilon}{3}<\varepsilon.
\end{aligned}
\]
This yields
\[
\sup_{x\in K}\|Tx-x\|<\varepsilon,
\quad \|T\|\leq\lambda,
\]
which implies that $X$ has $\lambda$-BAP.

We proceed with the proof of the converse direction by contradiction. Assume by contradiction that  $X$ does not have $\lambda$-BAP. Then there exists $(f_1,\ldots,f_m)\in X^m$ which does not belong to the norm closure of $\mathcal{C}_\lambda$. By the Hahn--Banach theorem, there are $g_1,\ldots,g_m\in X^*$ such that
\[
 \Rea\sum_jg_j(f_j)>
 \sup_{w=(w_1, \ldots, w_m)\in\mathcal C_\lambda}\Rea\sum_jg_j(w_j) = \sup_{T \in \F(X), \|T\| \leq \lambda}\Rea \ell(T),
\]
where $\ell$ is the finite test given below
$$
\ell(T)=\sum_{j=1}^m g_j(Tf_j),\quad  \forall\,  T \in \B(X). $$
% Over $\R$, $\Rea$ is omitted; over $\C$, every continuous real-linear separating functional is the real part of a continuous complex-linear functional.
 Let $a = \overline{\ell(T)}/|\ell(T)|$ if $\ell(T)\neq 0$, and when $\ell(T)=0$, $a = 0$ as convention. Then for any $T \in \F(X)$ with $\| T \| \leq \lambda$, one has
\[
|\ell(T)| = a \ell(T) = \ell(aT) = \Rea \ell(aT) \leq \sup_{S \in \F(X), \|S\| \leq \lambda}\Rea \ell(S).
\]
This implies
\begin{equation}\label{eq:2.1}
\sup_{S \in \F(X), \|S\| \leq \lambda} |\ell(S)| \leq \sup_{S \in \F(X), \|S\| \leq \lambda}\Rea \ell(S) < \Rea\sum_jg_j(f_j).
\end{equation}
However, by \eqref{h:equa2.1},
\[
\sup_{S \in \F(X), \|S\| \leq \lambda} |\ell(S)| = \lambda c_X(\ell) \geq |t_X(\ell)| \geq \Rea t_X(\ell) = \Rea\sum_jg_j(f_j),
\]
which leads to a contradiction with \eqref{eq:2.1}. This implies that $X$ has $\lambda$-BAP.

\iffalse
To pass to a compact set $K$ and error $\varepsilon>0$, choose a finite $\varepsilon/[3(\lambda+1)]$-net $\{f_j\}$ in $K$ and an admissible $T$ with error less than $\varepsilon/3$ on the net. If $x\in K$ and $\norm{x-f_j}<\varepsilon/[3(\lambda+1)]$, then
\[
 \norm{Tx-x}\leq(\lambda+1)\norm{x-f_j}+\norm{Tf_j-f_j}<2\varepsilon/3.
\]
This proves $\lambda$-BAP.
\fi

Now assume that $\beta(X)<\infty$. If $c_X(\ell)=0$ where $\ell $ is of the form \eqref{h:equa2.0}, take a bounded finite rank projection $P$ fixing $\spn\{f_j\}$. For $P\neq0$, apply the definition of $c_X(\ell)$ to $P/\norm{P}$ to obtain $\ell(P)=0$. Since $Pf_j=f_j$, this implies $t_X(\ell)=0$. In addition, the case $P=0$ is immediate. Hence, without loss of generality, we can assume $c_X(\ell) \neq 0$.
% A unit vector $f$ and a norm-one functional $g$ satisfying $g(f)=1$ give the test $\ell(T)=g(Tf)$. The rank-one contraction $Tx=g(x)f$ shows $c_E(\ell)=t_E(\ell)=1$.
It follows from \eqref{h:equa2.1} and the definition of $\beta(X)$ that $$\beta(X) \geq \sup_{c_X(\ell)>0; \ \ell \ \text{finite test}}\frac{\abs{t_X(\ell)}}{c_X(\ell)} = : \lambda_0. $$
Conversely, the first assertion shows that $X$ has $\lambda_0$-BAP, hence $\beta(X) \leq \lambda_0$, which proves \eqref{eq:test-beta} and that $X$ has $\beta(X)$-BAP.

The case $\beta(X) = \infty$ is immediate. We complete the proof.
\end{proof}

\subsection{Local transfer of finite rank dual operators} The following proposition will be used; see \cite[Theorem 1.1]{Oja2014}.

\begin{proposition}\label{h:prop2.2}
Let $M\subset X^{**}$ and $N\subset X^*$ be finite-dimensional subspaces, and let $\eta>0$. Then there exists a linear map $R:M\to X$ such that
\begin{equation}\label{eq:plr}
 \norm{R}\leq1+\eta,\quad \text{and} \quad  f(Ra)=a(f),
 \quad \forall \ a\in M,\ f\in N.
\end{equation}
\end{proposition}

The next lemma follows from Proposition \ref{h:prop2.2}.

\begin{lemma}\label{lem:local-transfer}
Let $K\in\F(X^*)$, $F_0\subset X^*$ be a finite-dimensional subspace, and $\eta>0$. Then there exists $U\in\F(X)$ such that for all $f\in F_0$,
\[
U^*f=Kf, \quad \norm{U}\leq(1+\eta)\norm{K}.
\]
% No weak-star continuity assumption is made on $K$.
\end{lemma}
\begin{proof}
If $K=0$, just take $U=0$. Then we suppose that $K\neq 0$. Choose a basis $y_1,\ldots,y_m$ of $\ran K\subset X^*$. Let $y_j^* \in X^{**}$ denote the coordinate functionals with respect to $y_j$. Hence
\[
 Kf=\sum_{j=1}^m y_j^*(Kf)y_j, \quad \forall\, f\in X^*.
\]
Denote by $a_j = y_j^* \circ K\in X^{**}$ for each $1\leq j\leq m$. Let $J_X$ denote the canonical embedding of $X$ into $X^{**}$. Then for each $x \in X$, one has
\[
 K^*J_Xx=\sum_{j=1}^m y_j(x)a_j\in M:=\spn\{a_1,\ldots,a_m\}\subset X^{**}.
\]
By using Proposition \ref{h:prop2.2} with $M$ and $N=F_0$, there exists a linear map $R:M\to X$ satisfying \eqref{eq:plr}. Let $U=RK^*J_X$. Then $U\in\F(X)$ and
\[
 \norm{U}\leq\norm{R}\norm{K^*}\norm{J_X}
 \leq(1+\eta)\norm{K}.
\]
By \eqref{eq:plr}, one has that for $x\in X$ and $f\in F_0$,
\[
 (U^*f)(x)=f(RK^*J_X x)=(K^*J_X x)(f)=(Kf)(x),
\]
as desired.
\end{proof}

\begin{lemma}\label{lem:interpolation}
Suppose that $X^*$ has $\lambda$-BAP for some $\lambda\geq1$. For every finite-dimensional subspace $F_0\subset X^*$ and $\varepsilon>0$, there exists $U\in\F(X)$ such that for each $f \in F_0$,
\[
U^*f=f, \quad \norm{U}\leq\lambda+\varepsilon.
\]
\end{lemma}
\begin{proof}
If $F_0=\{0\}$, just take $U=0$. Then we suppose that $F_0\neq 0$. Let $f_1,\ldots,f_m$ be a basis of $F_0$, and let $f_j^*$ denote the coordinate functionals with respect to $f_j$, extended to $X^{**}$. Choose $\eta >0$ such that $(1+\eta)(\lambda+\eta)\leq\lambda+\varepsilon$. Since $X^*$ has $\lambda$-BAP, there exists $H\in\F(X^*)$ of norm at most $\lambda$ such that
\[
 \sum_{j=1}^m\norm{f^*_j}\norm{f_j-Hf_j}\leq \eta.
\]
Therefore, the following finite rank operator
\[
 Tf=\sum_{j=1}^mf^*_j(f)(f_j-Hf_j)
\]
satisfies $\norm{T}\leq \eta$. Moreover, for every vector $f_j$, $Tf_j=f_j-Hf_j$. Hence, $K=H+T\in \F(X^*)$ fixes $F_0$ and has norm at most $\lambda+\eta$. Finally, the desired $U$ follows directly from Lemma~\ref{lem:local-transfer}.
\end{proof}
Lemma \ref{lem:interpolation} is classical and well-known to experts; it is also a special case of \cite[Corollary 4.1]{Oja2014}. We include the detailed proof just for the reader's convenience.
% Positivity will be obtained later from the graph norm and exact interpolation, not assumed of the dual approximations.
\iffalse
\subsection{Transport and compact-set approximation}
We record two elementary facts used below. First, if a uniformly bounded net of operators converges pointwise to the identity, then it converges uniformly on each norm-compact set. Indeed, a finite net in the compact set controls the pointwise error on finitely many vectors, and the common bound for $\norm{T_\alpha-\Id}$ controls the remaining error.

Second, if $W:E\to F$ is an isomorphism and $F$ has $\lambda$-BAP, then $E$ has $\lambda\norm{W}\norm{W^{-1}}$-BAP. For a compact $K\subset E$, approximate the identity on $W(K)$ by $T\in\F(F)$ and use $W^{-1}TW$. Its error on $K$ is at most $\norm{W^{-1}}$ times the error of $T$ on $W(K)$, and its norm is at most $\lambda\norm{W}\norm{W^{-1}}$.
\fi
\subsection{Property A}
In this subsection, let \(G\) be a finitely generated group with a finite symmetric generating set \(S\) which does not contain the identity of $G$. The left-invariant word metric on $G$ with respect to $S$ is defined as
\[
 d_G(x,y)=\abs{x^{-1}y}_S,
\]
where $|\cdot|_S$ denotes the word length with respect to $S$. For \(s\in G\), let \(L_s\delta_g=\delta_{sg}\) for any $g\in G$ on \(\ell_1(G)\). Denote by \(\Prob(G)\) the subset of all nonnegative vectors in \(\ell_1(G)\) with the $\ell_1$-norm \(1\).

Property A was first introduced by Yu \cite{Yu}. The following equivalent definition of Yu's Property A was proposed in \cite{HigsonRoe}. We refer the interested reader to \cite{NY23,Roe03} for a more detailed account of Property A.

\begin{definition}\label{def:a}
% In the probability-kernel formulation,
\(G\) has \emph{Property A} if for every \(R>0\) and \(\alpha>0\), there exist \(D<\infty\) and \(\xi:G\to\Prob(G)\) such that
\begin{itemize}
\item for any $x\in G$, $\supp\xi(x)\subset B(x,D);$
\item $d_G(x,y)\leq R\ \Longrightarrow\ \norm{\xi(x)-\xi(y)}_1\leq\alpha,$
\end{itemize}
where $B(x, D)$ denotes the ball in $G$ with center $x$ and radius $D$.
\end{definition}
%For word metrics, the balls are finite and their cardinality is independent of their center. Thus the bounded geometry hypothesis in the equivalence with the usual definition of Property A is satisfied \cite[Proposition 2.2]{Nowak}.

Let \(\mathcal M_{\mathrm{fin}}(G)\) denote the positive and total-summation-preserving operators \(A:\ell_1(G)\to\ell_1(G)\) whose columns $A\delta_g$ are all supported in one common finite subset of \(G\). Here, ``positive'' means that each entry of $A\delta_g$ is nonnegative for all $g\in G$. Every such operator has norm \(1\) and finite rank. Define
\begin{equation}\label{eq:kappa}
 \kappa_S(G)=
 \inf_{A\in\mathcal M_{\mathrm{fin}}(G)}
 \max_{s\in S}\norm{AL_s-L_sA}_{1\to1},
\end{equation}
where $\| \cdot \|_{1\to1}$ denotes the operator norm on $\ell_1(G)$. For the trivial group, we take the maximum in \eqref{eq:kappa} as \(0\) and  \(S=\varnothing\). As mentioned earlier in Introduction, $ \kappa_S(G)$ plays an essential role in the proofs of Theorem~\ref{thm:main-ap} and Theorem~\ref{thm:main-metric}.

The following proposition concerns an equivalent relationship between $ \kappa_S(G)$ and Property A.

\begin{proposition}\label{prop:markov}
 \(G\) has Property A if and only if \(\kappa_S(G)=0\). In particular, if \(G\) does not have Property A, then
\[
 0<\kappa_S(G)\leq2.
\]
\end{proposition}
\begin{proof}
The trivial group has Property A and $\kappa_S(G)=0$ by convention. Then we can assume that $S\neq\varnothing$. For a bounded operator $B:\ell_1(G)\to\ell_1(G)$, it is clear that
\begin{equation}\label{eq:column-norm}
 \norm{B}_{1\to1}=\sup_{g\in G}\norm{B\delta_g}_1.
\end{equation}
First, we show that \(\kappa_S(G)=0\) implies Property A of $G$. Let $A\in\mathcal M_{\mathrm{fin}}(G)$ and write $m_g=A\delta_g$. For $\delta>0$, suppose that
 \[ \max_{s \in S} \| AL_s - L_sA \|_{1 \to 1} \leq \delta. \]
 Let $J$ be a common finite support for all $m_g$. Define
\begin{equation}\label{h:xi}
 \xi(x)=L_xm_{x^{-1}},\quad D=\max_{j\in J}|j|_S<\infty.
\end{equation}
% Translations preserve positivity and mass.
Then $\xi(x)\in \Prob(G)$. Since $d_G(x,xj)=|j|_S$, one has
\[
 \supp\xi(x)\subset xJ\subset B(x,D).
\]
Since $L_{x}$ is an isometry, for any $s\in S$
\begin{align}\label{eq:2.5}
 \norm{\xi(xs)-\xi(x)}_1
 &=\norm{L_xL_sm_{s^{-1}x^{-1}}-L_x m_{x^{-1}}}_1 \notag\\
 &=\norm{(L_sA-AL_s)\delta_{s^{-1}x^{-1}}}_1\leq\delta.
\end{align}
Note that for every $x, y \in G$ with $d_G(x,y)\leq R$, one can write $x^{-1}y$ as a word in $S$ of length at most $\lceil R\rceil$, namely
\[
x^{-1}y = \prod_{j} s_j, \quad s_j \in S,
\]
where the number of $s_j$ appearing on the right-hand-side is at most $\lceil R\rceil$. Then by \eqref{eq:2.5}, we obtain that $d_G(x,y)\leq R$ implies
\[
\| \xi(x) - \xi(y) \|_1 = \| \xi(x) - \xi(x (x^{-1}y) ) \|_1 \leq \lceil R\rceil \delta.
\]
% If $d_G(x,y)\leq R$, write $x^{-1}y$ as a word in $S$ of length at most $\lceil R\rceil$ and sum these edge estimates.
If $\kappa_S(G)=0$, for any given $R>0$ and $\alpha>0$, one can choose an $A$ with $0<\delta \leq \alpha/\lceil R\rceil$. Hence, the resulting kernel $\xi$ in \eqref{h:xi} satisfies Definition~\ref{def:a}, and this implies Property A of $G$.

Conversely, suppose that $G$ has property A. For $\delta>0$, take a kernel $\xi$ such that $\exists\, D>0$,
\[
\supp\, \xi(x) \subset B(x, D) \quad \text{and} \quad \| \xi(x) - \xi(y) \|_1 \leq \delta, \quad \forall\, d_G(x, y) \leq 1.
\]
For each $g \in G$, let $m_g=L_g\xi(g^{-1})$ and define $A\in\mathcal M_{\mathrm{fin}}(G)$ via
\[
  \quad A\delta_g=m_g.
\]
Since $gB(g^{-1},D)=B(e,D)$, all columns $A\delta_g$ have support in the same finite ball $B(e,D)$. Then it is clear that $A$ extends to an operator in $\mathcal{M}_{\rm fin}(G)$. For $s\in S$, since $L_{sg}$ is an isometry and $d_G(g^{-1}s^{-1},g^{-1})=|s|_S\leq1$,
\begin{align*}
 \norm{m_{sg}-L_sm_g}_1 &= \norm{ L_{sg}\xi(g^{-1}s^{-1}) - L_{sg}\xi(g^{-1}) }_1\\
 &=\norm{\xi(g^{-1}s^{-1})-\xi(g^{-1})}_1\leq\delta.
\end{align*}
Therefore, \eqref{eq:column-norm} yields $\norm{AL_s-L_sA}_{1\to1}\leq\delta$. By letting $\delta\downarrow0$, one can conclude that $\kappa_S(G)=0$.

Finally, it is clear that $\kappa_S(G)\leq 2$. This completes the proof.
\end{proof}

%The columns \(m_g\) depend on \(g\), as do the kernels in the definition of Property A. This dependence is essential; replacing them by one common probability vector would give an amenability condition.

The values of \(\kappa_S(G)\) may vary from the generating sets. For a nontrivial group, if \(S'\) is another finite symmetric generating set and \(L=\max_{s\in S}\abs{s}_{S'}\), expansion of the commutator along a word gives \(\kappa_S(G)\leq L\kappa_{S'}(G)\). Reversing the roles of \(S\) and \(S'\) shows that $\kappa_S(G)=0$ is independent of this choice of generating sets, and so is Property A of $G$.

\bigskip

\section{A dual space with BAP but without MAP}\label{sec:construction}

In this section, we construct a dual space of a MAP Banach space, which has BAP but fails MAP. The construction relies on a finitely generated group without property A.

\subsection{Main construction}
% There are two independent stages: construction of a dual with $1<\beta<\infty$, and amplification within the class of duals.
Let $G$ be a finitely generated group and $S$ be a finite symmetric generating set that does not contain the identity of $G$. Now we encode the generators of \(G\) via a bounded operator (denoted by $Q$ in the following) on labelled copies of \(\ell_1(G)\).
% The labels will keep the probability columns separate while the operator records their relations.
We use distinct formal labels of $t_s$ for $s\in S$, and define
\[
 \mathcal T=\{0,1,d\}\cup\{t_s:s\in S\},\quad I=\mathcal T\times G.
\]
The labels \(0, 1, d\) are separate from all labels \(t_s\). Let
\[
 Z=\ell_1(G)\oplus_1\ell_1(G).
\]
Define \(Q:\ell_1(I)\to Z\) on the standard basis $(e_i)_{i\in I}\subset \ell_1(I)$ by
\begin{equation}\label{eq:Q}
\begin{aligned}
 Qe_{0,g}& =(\delta_g,0), \quad
 Qe_{1,g} =(0,\delta_g),\\
 Qe_{d,g}&=\frac12(\delta_g,\delta_g), \quad
 Qe_{t_s,g} =\frac12(\delta_g,L_s\delta_g).
\end{aligned}
\end{equation}
Every column is nonnegative and has \(Z\)-norm \(1\). Moreover, \(\norm{Q}_{\ell_1(I) \to Z}=1\). Let \(\sigma_I\) and \(\sigma_Z\) denote the canonical summation functionals on $\ell_1(I)$ and $Z$ respectively. Then, we have
\begin{equation}\label{eq:mass}
 \sigma_ZQ=\sigma_I.
\end{equation}

Fix \(0<a\leq1\). Define
\begin{equation*}\label{eq:3.2}
X_{G,a}=(\ell_1(I),n_a),\quad \text{and} \quad
 n_a(x)=\norm{x}_1+a\norm{Qx}_Z,\quad \forall\, x\in \ell_1(I).
\end{equation*}
% where $Z=\ell_1(G)\oplus_1\ell_1(G)$ and the positive operator $Q$ encodes the generators.
% Coordinate projections are contractions for $n_a$, and
Then for any $(j, g) \in I$ and $x \in X_{G, a}$,
\[
n_a(e_{j, g}) = 1 + a, \quad \norm{x}_1\leq n_a(x)\leq(1+a)\norm{x}_1.
\]
This implies that $X_{G,a}$ is isomorphic to $\ell_1(I)$.
\iffalse
Thus $X_{G,a}$ has MAP and $\beta(E_{G,a})\leq1+a$.

Write $L_s\delta_g=\delta_{sg}$ and set
\[
 \kappa_S(G)=\inf_{A\in\mathcal M_{\mathrm{fin}}(G)}
 \max_{s\in S}\norm{AL_s-L_sA}_{1\to1},
\]
where $\mathcal M_{\mathrm{fin}}(G)$ consists of positive mass-preserving operators whose columns share a finite output support. We prove that $\kappa_S(G)=0$ exactly when $G$ has Property A. Local reflexivity transfers finite-dimensional interpolation by arbitrary finite-rank dual operators to the predual. Exact interpolation of the label indicators fixes the signed mass of each column on each label.

An operator with norm close to $1$ then has columns close to probabilities on the corresponding labels. Finite rank permits a uniform truncation, and relations in $\ker Q$ produce the translation-commutator estimate. The resulting bounds are
\begin{equation}\label{eq:intro-ca}
 C_a=1+(1+a)^2,\qquad b_a=\frac{a\kappa_S(G)}{8C_a},
\end{equation}
\begin{equation}\label{eq:intro-bounds}
 1+b_a\leq\beta(E_{G,a})\leq1+a.
\end{equation}
Osajda's theorem provides finitely generated groups without Property A \cite[Theorem 2]{Osajda}. For such a group $b_a>0$, giving Theorem~\ref{thm:main-metric}.
\fi
Its canonical dual, denoted by $E_{G,a}=X_{G,a}^*$, is
\begin{equation*}\label{eq:q_a}
 E_{G,a}=(\ell_\infty(I),q_a),\quad
 q_a(f):=\sup_{n_a(x)\leq1}\abs{\sum_{i\in I}f(i)x(i)},
\end{equation*}
and
\begin{equation}\label{eq:q_a-equivalent}
 \frac{1}{1+a}\norm{f}_\infty\leq q_a(f)\leq\norm{f}_\infty.
\end{equation}
All operator norms on \(E_{G,a}\) will be taken with respect to \(q_a\). The norm equivalence implies the BAP of $E_{G, a}$.
% while the form of \(n_a\) preserves contractive coordinate projections on the predual.

\begin{proposition}\label{prop:basic}
For every finitely generated \(G\) and \(0<a\leq1\), the space \(E_{G,a}\) has \((1+a)\)-BAP and the space \(X_{G,a}\) has MAP, over either scalar field.
\end{proposition}
\begin{proof}
We first prove that $E_{G,a}$ has the $(1+a)$-BAP. Let $f_1,\ldots,f_m\in E_{G,a}=\ell_\infty(I)$ and $\varepsilon>0$. Consider the map
	\[
	\Phi:I\longrightarrow\mathbb K^m,\qquad
	\Phi(i)=(f_1(i),\ldots,f_m(i)),
	\]
where $\mathbb K^m$ is equipped with the maximum norm
	\[
	\|(z_1,\ldots,z_m)\|_\infty
	=\max_{1\le j\le m}|z_j|.
	\]
Since each $f_j$ is bounded, the set
	\[
	A=\Phi(I)
	=\{(f_1(i),\ldots,f_m(i)):i\in I\}
	\]
is bounded in the finite-dimensional space $\mathbb K^m$. Consequently, $A$ is relatively compact. Thus, there exist
	$v_1,\ldots,v_N\in\mathbb K^m$ such that
	\[
	A\subseteq
	\bigcup_{k=1}^N B_\infty(v_k,\varepsilon/3),
	\]
where
	\[
	B_\infty(v,r)
	=\{z\in\mathbb K^m:\|z-v\|_\infty<r\}.
	\]
For $1\le k\le N$, define
	\[
	E_k=\Phi^{-1}\bigl(B_\infty(v_k,\varepsilon/3)\bigr).
	\]
Then $I=\bigcup_{k=1}^N E_k$. To obtain a disjoint partition, define via the following recursive law
	\[
	C_1=E_1,\qquad
	C_k=E_k\setminus\bigcup_{l=1}^{k-1}E_l,
	\quad 2\le k\le N.
	\]
Discarding the empty sets, we obtain a finite partition
	\[
	\mathcal P=\{C_1',\ldots,C_M'\}
	\]
of $I$ into a pairwise disjoint nonempty union. For each $C\in\mathcal P$, choose a representative $i_C\in C$. Since every $C$ is contained in some $E_k$, for any $i\in C$ we have
	\[
	\begin{aligned}
		\|\Phi(i)-\Phi(i_C)\|_\infty
		&\le \|\Phi(i)-v_k\|_\infty
		+\|v_k-\Phi(i_C)\|_\infty\\
		&<\frac{2\varepsilon}{3}
		<\varepsilon.
	\end{aligned}
	\]
Thus, for every $C\in\mathcal P$,
	\[
	\sup_{i\in C}
	\max_{1\le j\le m}|f_j(i)-f_j(i_C)|
	<\varepsilon.
	\]
Define $P:\ell_\infty(I)\to\ell_\infty(I)$ by
	\[
	(Pf)(i)=f(i_C),\qquad i\in C,\quad C\in\mathcal P.
	\]
It is clear that $P \in \F(\ell_\infty(I))$ and
\[
 \norm{P}_{\infty\to\infty}\leq1,\quad
 \norm{Pf_j-f_j}_\infty<\varepsilon.
\]
By \eqref{eq:q_a-equivalent}, one has
\[
 q_a(Pf)\leq\norm{f}_\infty\leq(1+a)q_a(f),\quad
 q_a(Pf_j-f_j)<\varepsilon.
\]
By the same argument as in Lemma~\ref{lem:finite-tests}, we obtain that $E_{G,a}$ has the $(1+a)$-BAP. In particular, when $a=0$, this also proves the usual MAP of $\ell_\infty(I)$.

For finite $J\subset I$, let
$P_J:\ell_1(I)\to\ell_1(I)$ denote the canonical
coordinate projection defined by
\[
(P_Jx)(i)=
\begin{cases}
	x(i), & i\in J,\\
	0, & i\notin J.
\end{cases}
\]
Then $P_J\in \F(X_{G,a})$. By using $\norm{Q}_{\ell_1(I)\to Z}=1$ and the disjoint supports of $P_Jx$ and $(\Id-P_J)x$, we obtain
\begin{align*}
 n_a(P_Jx)
 &\leq\norm{P_Jx}_1+a\norm{Qx}_Z+a\norm{Q(\Id-P_J)x}_Z\\
 &\leq\norm{P_Jx}_1+a\norm{Qx}_Z+a\norm{(\Id-P_J)x}_1\\
 &=n_a(x)-(1-a)\norm{(\Id-P_J)x}_1\leq n_a(x).
\end{align*}
As $J$ increases through the finite subsets of $I$,
\[
 n_a(x-P_Jx)\leq(1+a)\norm{x-P_Jx}_1\longrightarrow0.
\]
Hence, $(P_J)$ is a net of finite rank contractions converging pointwise, and therefore uniformly on norm-compact sets to the identity on $X_{G,a}$. This yields that $X_{G, a}$ has MAP.
% Its adjoints converge pointwise in the weak-star topology on $E_{G,a}$; this alone does not assert norm approximation on compact subsets of the dual.
\end{proof}

The following theorem concerns the optimal BAP constant of $E_{G,a}$ to the obstruction $\kappa_S(G)$, which is the key ingredient in the proof of Theorem \ref{thm:main-metric}.
\begin{theorem}\label{thm:3.2}
Let $E_{G, a}$ be defined as above, and
\[
b_a=\frac{a\kappa_S(G)}{8(1 + (1 + a)^2)}.
\]
Then
\begin{equation}\label{eq:intro-bounds}
 1+b_a\leq\beta(E_{G,a})\leq1+a.
\end{equation}
\end{theorem}

The proof of Theorem \ref{thm:3.2} is quite technical, and we put the detailed proof to the end of this section. Now we come to the proof of Theorem \ref{thm:main-metric}.

\subsection{Proof of Theorem \ref{thm:main-metric}}
By virtue of Theorem \ref{thm:3.2}, we show directly Theorem~\ref{thm:main-metric}.

\begin{proof}[Proof of Theorem~\ref{thm:main-metric}]
By \cite{Gro03,Osajda,Osa18}, one can choose a finitely generated
group $G$ without Property A.
Then by Proposition \ref{prop:markov}, $\kappa_S(G)>0$, and by Theorem \ref{thm:3.2}, we obtain
\[
 1<\beta(E_{G,a})\leq1+a.
\]
Together with the previously established properties of $X_{G,a}$ in Proposition \ref{prop:basic}, this gives the desired conclusion.
\end{proof}

\subsection{Proof of Theorem \ref{thm:3.2}}
For \(t\in\mathcal T\), let \(f_t\) be the indicator function of \(\{t\}\times G\) in \(E_{G,a}\), and denote by
\[
 F_0=\spn\{f_t:t\in\mathcal T\},\quad \mathbf1=\sum_{t\in\mathcal T}f_t.
\]
We consider finite rank operators $U$ on $X_{G,a}$ satisfying $U^*f_t=f_t$ for every $f_t$ as \(t\in\mathcal T\). These identities require $U$ to preserve the total summation on each label \(t\in\mathcal T\). Throughout this section, \(C_a=1+(1+a)^2\).

We say an operator $A$ on $\ell_1(I)$ is label-preserving if for any $i=(t,g)\in I$, the support of $Ae_i$ is contained in $\{t\}\times G$. We need the following two auxiliary lemmas below.

\begin{lemma}\label{lem:extraction}
Over the scalar field \(\R\), let \(\delta\geq0\) and suppose that \(U\in\F(X_{G,a})\) satisfies
\[
 \norm{U}_{n_a\to n_a}\leq1+\delta,\quad U^*f_t=f_t,\quad  t\in\mathcal T,
\]
where $\| \cdot \|_{n_a\to n_a}$ denotes the operator norm on $X_{G,a}$. Then, for every \(0<\theta<1\), there exists a positive, total-summation-preserving, label-preserving operator \(S_0\) on \(\ell_1(I)\), with a common finite support of all $S_0e_i$ for $i\in I$, such that
\begin{equation}
 \norm{S_0-U}_{1\to1} \leq(1+a)\delta+2\theta,\label{eq:column-error}
\end{equation}
and
\begin{equation}
\norm{S_0}_{n_a\to n_a} \leq1+C_a\delta+2(1+a)\theta.\label{eq:S-estimate}
\end{equation}
\end{lemma}
\begin{proof}
Denote by $u_i=Ue_i$ with $i=(t,g)\in I$. Since $U^*f_t=f_t$ for every $t\in\mathcal T$, we have
\begin{equation}\label{eq:label-mass}
 \sum_{h\in G}u_{t,g}( r,h )=\delta_{r,t},
 \quad r,t\in\mathcal T,\ g\in G.
\end{equation}
Consequently, together with \eqref{eq:mass},
$$\sigma_I(u_i)=1, \quad \norm{Qu_i}_Z\geq|\sigma_Z(Qu_i)|=|\sigma_I(u_i)|=1.$$
Since $ \norm{U}_{n_a\to n_a}\leq1+\delta$ and $n_a(e_i)=1+a$,
\[
 \norm{u_i}_1+a\leq n_a(u_i) \leq (1 + \delta) n_a(e_i) = (1+\delta)(1+a).
\]
This implies
\begin{equation}\label{eq:a-size}
 \norm{u_i}_1\leq1+(1+a)\delta, \quad i\in I.
\end{equation}
Hence, $\{ u_i : i \in I \}$ is bounded in the finite-dimensional subspace $\ran U$, and so its closure is compact in $\ell_1(I)$. Denote by $\Gamma_{\theta/2}$ a finite $(\theta/2)$-net for this closure in $\ell_1(I)$. For each $f \in \Gamma_{\theta/2}$, define the coordinate $I_f \subset I$ as a finite subset such that
\[
\sum_{i \in I\setminus I_f} |f(i)| < \frac{\theta}{2}.
\]
Let $J = \cup_{f \in \Gamma_{\theta/2}} I_f $. Since $\norm{\Id-P_J}_{1\to1}\leq1$, one has
\begin{align}\label{eq:tail}
 \sup_i\norm{(\Id-P_J)u_i}_1&\leq  \sup_i\norm{(\Id-P_J)(u_i-f_{u_i})}_1+ \sup_i\norm{(\Id-P_J)f_{u_i}}_1 \notag\\
 &\leq \sup_i\norm{u_i-f_{u_i}}_1+ \frac{\theta}{2}\leq \theta,
\end{align}
where $u_i$ belongs to the ball with center $f_{u_i}\in \Gamma_{\theta/2}$ and radius $\theta/2$ in $\ell_1(I)$. The finite set $J$ may depend on $U$ and $\theta$. Fix $i=(t,g)$, define $b_i \in \ell_1(I)$ by
\[
b_i (j) = \begin{cases}
      u_i (j), \quad & j \in J \cap (\{t\}\times G) \text{ and } u_i (j)> 0; \\
      0, & \text{ otherwise}.
\end{cases}
\]
Then by \eqref{eq:label-mass} and \eqref{eq:tail},
\begin{equation}\label{h:eq:bi}
\| u_i \|_1 \geq \| b_i \|_1 \geq 1 - \| (\Id - P_J)u_i \|_1 \geq 1 - \theta > 0.
\end{equation}
Therefore, the vector $p_i=b_i/ \| b_i \|_1$ is a probability supported on its own label in $J$. The definition of $b_i$ and $p_i$ requires real scalars $\mathbb{R}$.

Note that
\[
 \norm{u_i-b_i}_1=\norm{u_i}_1- \|b_i\|_1,
 \quad \norm{b_i-p_i}_1=| \| b_i \|_1 -1|.
\]
If $\| b_i \|_1 \leq1$, by \eqref{eq:a-size} and \eqref{h:eq:bi}, one has
\[
\| u_i - p_i \|_1 \leq \norm{u_i}_1- \|b_i\|_1 + 1 - \| b_i \|_1 \leq (1 + a)\delta + 2\theta.
\]
If $\| b_i \|_1 \geq 1$, by \eqref{eq:a-size}, one has
\[
\| u_i - p_i \|_1 \leq \| u_i \|_1 - 1 \leq (1 + a)\delta.
\]
Hence in both cases
\[
 \norm{u_i-p_i}_1\leq(1+a)\delta+2\theta.
\]
Define $$ S_0x=\sum_{i \in I} x(i) p_i, \quad x\in\ell_1(I).$$
This series is absolutely convergent, since $\norm{p_i}_1=1$. It is clear that $S_0$ is linear, positive, total-summation-preserving and label-preserving, and with a common finite support $J$ of all $S_0e_i$. In addition, for each $i\in I$, one has
\[
\| S_0e_i - Ue_i \|_1 = \| p_i - u_i \|_1 \leq (1+a)\delta+2\theta.
\]
By \eqref{eq:column-norm}, this implies \eqref{eq:column-error}. Moreover, for every $x$,
\begin{align*}
 n_a((S_0-U)x)
 &\leq(1+a)\norm{(S_0-U)x}_1\\
 &\leq(1+a)((1+a)\delta+2\theta)\norm{x}_1\\
 &\leq(1+a)((1+a)\delta+2\theta)n_a(x).
\end{align*}
Combining this with $n_a(Ux)\leq(1+\delta)n_a(x)$, this completes the proof of \eqref{eq:S-estimate}.
\end{proof}

%The relations in \(\ker Q\) now compare the probability columns belonging to different labels. Their disjoint supports preserve the \(\ell_1\) norm of each relation vector, so any increase in the graph norm must come from its image under \(Q\).

The following lemma concerns the upper bound of $\kappa_S(G)$.

\begin{lemma}\label{lem:commutators}
Under the same assumptions as in Lemma~\ref{lem:extraction}, for every \(0<\theta<1\), there exists \(A\in\mathcal M_{\mathrm{fin}}(G)\) such that
\begin{equation}\label{eq:commutator-bound}
 \max_{s\in S}\norm{AL_s-L_sA}_{1\to1}
 \leq\frac{8C_a}{a}\delta+\frac{16(1+a)}{a}\theta.
\end{equation}
\end{lemma}
\begin{proof}
It is immediate when the generating set $S=\varnothing$. Now suppose that $S\neq \varnothing$. Let $S_0$ be given as in Lemma~\ref{lem:extraction}, and $\omega=C_a\delta+2(1+a)\theta$. For  \(t\in\mathcal T\), let $\iota_t:\ell_1(G)\to\ell_1(I)$ denote the isometric inclusion $\iota_t : \delta_g \mapsto e_{t, g}$. By Lemma~\ref{lem:extraction}, there are probability vectors $p_{t,g}\in \Prob(G)$ with a common finite support such that
\[
 S_0e_{t,g}=\iota_t p_{t,g}.
\]
Define $A\delta_g=p_{0,g}$ and $D\delta_g=p_{1,g}$. Both $A$ and $D$ belong to $\mathcal M_{\mathrm{fin}}(G)$.

For $t=d$, let $s=e$ the identity of $G$, and for $t=t_s$ use the corresponding generator $s$. Here, we do not need to define the associated $s$ for $t=0$ or $1$. For $t\in \{d, t_s: s\in S\}$, define
\[
 h_{t,g}=e_{t,g}-\frac12e_{0,g}-\frac12e_{1,sg}.
\]
Then \eqref{eq:Q} implies $Qh_{t,g}=0$. Note that $n_a(h_{t,g})=\norm{h_{t,g}}_1=2$. In addition, the images of $e_{t,g}, e_{0,g}, e_{1,sg}$ under $S_0$ have disjoint label supports, and each column is a probability. Consequently,
\[
 \norm{S_0h_{t,g}}_1=1+\frac12+\frac12=2.
\]
By \eqref{eq:S-estimate}, we have
\[
 2+a\norm{QS_0h_{t,g}}_Z=n_a(S_0h_{t,g})\leq (1+\omega)n_a(h_{t,g}) \leq 2(1 + \omega).
\]
This implies
\[
\norm{QS_0h_{t,g}}_Z\leq \frac{2\omega}{a}.
\]
For $t= t_s$, by definition of $Q$, one has
\begin{equation*}\label{eq:kernel-image}
 QS_0h_{t,g}=\frac12(p_{t, g}-A\delta_g,\ L_sp_{t, g}-D\delta_{sg}).
\end{equation*}
This yields
\begin{equation}\label{h:eq:ps}
 \norm{p_{t, g}-A\delta_g}_1+\norm{L_sp_{t, g}-D\delta_{sg}}_1\leq \frac{4\omega}{a}.
\end{equation}
Since $L_s$ is an isometry,
\[
 \norm{D\delta_{sg}-L_sA\delta_g}_1 \leq \norm{L_sp_{t, g}-L_sA\delta_g}_1 +  \norm{L_sp_{t, g}-D\delta_{sg}}_1 \leq \frac{4\omega}{a}.
\]
This implies
\begin{equation}\label{eq:3.12}
\| DL_s - L_sA \|_{1 \to 1} \leq \frac{4\omega}{a}.
\end{equation}
For $t=d$, in \eqref{h:eq:ps}, one has $s=e$ and the triangle inequality yields $ \norm{A\delta_g-D\delta_g}_1\leq4\omega/a $. This implies \begin{equation}\label{eq:3.13}
\norm{A-D}_{1\to1}\leq \frac{4\omega}{a}.
\end{equation}
Then by \eqref{eq:3.12} and \eqref{eq:3.13},
\[
 \norm{AL_s-L_sA}_{1 \to 1}
 \leq\norm{(A-D)L_s}_{1 \to 1}+\norm{DL_s-L_sA}_{1 \to 1} \leq \frac{8\omega}{a}.
\]
This finishes the proof of \eqref{eq:commutator-bound}.
\end{proof}

Recall that $F_0=\spn\{f_t:t\in\mathcal T\}$. Now we can give a lower bound for the operator norm of $U$ on $X_{G,a}$ whose dual operator fixes $F_0$ in terms of $\kappa_S(G)$.

\begin{corollary}\label{cor:exact-obstruction}
Over either scalar field, every \(U\in\F(X_{G,a})\) satisfying \(U^*|_{F_0}=\Id\) satisfies
\begin{equation}\label{eq:exact-lower}
 \norm{U}_{n_a\to n_a}\geq1+b_a,\quad b_a=\frac{a\kappa_S(G)}{8C_a}.
\end{equation}
\end{corollary}
\begin{proof}
The nonzero constant functional $\mathbf1=\sum_tf_t$ is fixed by $U^*$, so $$\norm{U}_{n_a\to n_a}=\norm{U^*}_{q_a\to q_a}\geq1.$$
Over $\R$, let $\delta=\norm{U}_{n_a\to n_a}-1\geq 0$. For every $0<\theta<1$, by Lemma~\ref{lem:commutators}, there exists an operator $A \in \mathcal{M}_{\rm fin}(G)$ satisfying \eqref{eq:commutator-bound} and therefore
\[
 \kappa_S(G)\leq\frac{8C_a}{a}\delta+\frac{16(1+a)}a\theta.
\]
% This inequality holds for every $\theta>0$, although the output supports may vary with $\theta$.
By letting $\theta\downarrow0$, we obtain $\delta\geq a\kappa_S(G)/(8C_a) = b_a$, and this proves \eqref{eq:exact-lower} over $\R$. For the complex scalar field $\C$, let $X_{G,a}^{\R}$ denote the real coordinate subspace with the same graph norm and define the real-linear operator by
\[
 U_{\R}:X_{G,a}^{\R}\to X_{G,a}^{\R},\quad U_{\R}x=\Rea(Ux).
\]
It is clear that $U_{\R}$ has real finite rank if $U$ is of finite rank.  Since $Q$ has real coefficients, $Q\Rea z=\Rea Qz$ and
\[
 n_a(\Rea z)\leq n_a(z).
\]
This yields $\norm{U_{\R}}_{n_a\to n_a}\leq\norm{U}_{n_a\to n_a}$. For every real $x$,
\[
 f_t(U_{\R}x)=\Rea f_t(Ux)=\Rea f_t(x)=f_t(x).
\]
Thus $U_{\R}$ satisfies the hypotheses of Lemma \ref{lem:extraction} and Lemma \ref{lem:commutators} with $\delta=\norm{U}_{n_a\to n_a}-1\geq0$. The same argument yields the desired required lower bound of $\norm{U}_{n_a\to n_a}$.
\end{proof}

\iffalse
We apply local reflexivity to transfer the preceding obstruction to the dual. A finite-rank operator that nearly fixes the label indicators can be corrected to one that fixes their span exactly. The correction also gives a quantitative approximation error.
\fi
We are ready to show Theorem \ref{thm:3.2}.

\begin{proof}[Proof of Theorem \ref{thm:3.2}]
By Proposition~\ref{prop:basic},
$$ \beta(E_{G,a})\leq1+a.$$
Suppose that $E_{G,a}$ has $\lambda$-BAP. By Lemma~\ref{lem:interpolation}, for every $\eta>0$, there exists a finite rank operator $U$ on $X_{G,a}$ such that
$$U^*|_{F_0}=\Id, \quad \norm{U}_{n_a\to n_a}\leq\lambda+\eta .$$
By Corollary~\ref{cor:exact-obstruction}, we obtain $1+b_a\leq\lambda+\eta$. By letting $\eta\downarrow0$ and then taking the infimum over admissible $\lambda$, we have
$$   1+b_a\leq\beta(E_{G,a}). $$
This completes the proof.
\end{proof}

Moreover, we can obtain a quantitative estimate for failing the approximation by finite rank operators.
\begin{corollary}
Let \(\lambda\geq1\) and \(k=\abs{\mathcal T}\). If \(T\in\F(E_{G,a})\) and \(\norm{T}\leq\lambda\), then
\begin{equation}\label{eq:finite-gap}
 \max_{t\in\mathcal T}q_a(Tf_t-f_t)
 \geq\frac{ \max\{1+b_a-\lambda, 0 \} }{(1+a)k}.
\end{equation}
\end{corollary}
\begin{proof}
The $f_t$'s are linearly independent and
\[
 \norm{\sum_r\alpha_rf_r}_\infty=\max_r|\alpha_r|.
\]
Thus the coordinate functional $\phi_t(\sum_r\alpha_rf_r)=\alpha_t$ on $F_0$ has norm at most $1+a$ with respect to $q_a$ by \eqref{eq:q_a-equivalent}. By the Hahn--Banach theorem, $\phi_t$ extends to a functional in $E_{G,a}^*$. Denote by
\[
 \Delta=\max_{t} q_a(Tf_t-f_t),\quad
 K=T+\sum_t\phi_t(\cdot)(f_t-Tf_t).
\]
Then $K \in \F(E_{G, a})$ and $Kf_r=f_r$, hence $K|_{F_0}=\Id$. Moreover,
\[
 \norm{K}\leq\norm{T}+\sum_t\norm{\phi_t}q_a(f_t-Tf_t)
 \leq\lambda+(1+a)k\Delta.
\]
For every $\delta > 0$, by Lemma~\ref{lem:local-transfer}, there exists $U\in\F(X_{G,a})$ with $U^*|_{F_0}=\Id$ and $\norm{U}_{n_a\to n_a}\leq(1+\delta)\norm{K}$. By Corollary~\ref{cor:exact-obstruction}, one has
\[
 1+b_a\leq(1+\delta)\bigl(\lambda+(1+a)k\Delta\bigr).
\]
This finishes the proof by letting $\delta\downarrow 0$.
\end{proof}

In particular, a contractive approximation on the dual would force the group to have Property A.
\begin{corollary}\label{thm:group}
If \(E_{G,a}\) has MAP, then \(G\) has Property A. If \(G\) fails Property A, every finite rank contraction on \(E_{G,a}\) has error at least \(b_a/((1+a)\abs{\mathcal T})>0\) on the label indicators $f_t$.
\end{corollary}
\begin{proof}
MAP of $E_{G, a}$ implies \(\beta(E_{G,a})=1\). Then \eqref{eq:intro-bounds} forces \(\kappa_S(G)=0\), and Proposition~\ref{prop:markov} implies Property A of $G$. The second assertion follows from \eqref{eq:finite-gap} with \(\lambda=1\).
\end{proof}

\bigskip

\section{A dual with AP and without BAP}\label{sec:tensor}
The preceding construction gives \(1<\beta(E_{G,a})<\infty\). To construct dual Banach spaces with AP but without BAP, we need dual spaces with arbitrarily large finite constants $\lambda$. The desired example can be obtained by projective tensor product of Banach spaces. We proceed with Theorem~\ref{thm:main-tensor}.

\subsection{Proof of Theorem~\ref{thm:main-tensor}}
The following theorem gives a quantitative lower bound for the optimal
bounded approximation constants of duals of projective tensor products.
It plays a central role in deriving Theorem~\ref{thm:main-ap}
from Theorem~\ref{thm:main-metric}.
\begin{theorem}\label{thm:main-tensor}
	For nonzero Banach spaces $X$ and $Y$ over the same scalar field,
	\begin{equation}\label{eq:intro-tensor}
		\beta\bigl((X\tens Y)^*\bigr)\geq\beta(X^*)\beta(Y^*).
	\end{equation}
	The inequality holds for extended values in $[1,\infty]$.
\end{theorem}

Theorem \ref{thm:main-tensor} allows a strict gap above $1$ to be amplified
exponentially by taking tensor powers. Applied to the space in
Theorem~\ref{thm:main-metric}, it yields dual spaces with finite
but unbounded BAP constants.

% Norm-one slice embeddings and projections give the elementary bound by $\max\{\beta(X^*),\beta(Y^*)\}$. Theorem~\ref{thm:main-tensor} replaces that maximum by a product.
Since $X^*$ and $Y^*$ are isometric to norm-one complemented
subspaces of $(X\widehat\otimes_\pi Y)^*$ respectively, one immediately obtains
\[
\beta\bigl((X\widehat\otimes_\pi Y)^*\bigr)
\geq \max\{\beta(X^*),\beta(Y^*)\}.
\]
Our theorem strengthens this elementary maximum bound
to a multiplicative lower bound.
% Note that it is separate from the familiar upper estimate for approximation constants on projective tensor products themselves; see \cite[Exercise 4.5]{Ryan} and its quantitative results.

We turn to the proof of Theorem~\ref{thm:main-tensor}. Let \(X,Y\neq\{0\}\) and denote by
\[
 D=(X\tens Y)^*
\]
with the bounded bilinear forms on \(X\times Y\), using the standard isometric duality for the projective tensor norm \cite[Chapter 2]{Ryan}.
For \(u\in Y^*\) and \(y\in Y\), define
\begin{align*}
 J_u:X^* \to D,&\quad (J_uf)(x,y')=f(x)u(y');\\
 S_y:D \to X^*,&\quad (S_yh)(x)=h(x,y).
\end{align*}
Then \(\norm{J_u}\leq\norm{u}\) and \(\norm{S_y}\leq\norm{y}\). For \(g\in X^{**}\), define
\begin{equation}\label{eq:slice-T}
 T_g:D \to Y^*,\quad (T_gh)(y)=g(S_yh).
\end{equation}
The right-hand-side of \eqref{eq:slice-T} is a linear functional on \(Y\), with absolute value at most \(\norm{g}\norm{h}\norm{y}\). Thus \(T_g\) takes values in \(Y^*\) and \(\norm{T_g}\leq\norm{g}\) for every \(g\in X^{**}\). The slice maps allow us to transfer a finite rank operator on \(D\) to one on \(Y^*\), which is exactly the following lemma.

\begin{lemma}\label{lem:compression}
Let \(\ell_{X^*}(A)=\sum_{j = 1}^m g_j(Af_j)\) be a finite test on \(\B(X^*)\) where $f_j\in X^*$ and $g_j\in X^{**}$. For \(H\in\B(D)\), define
\begin{equation*}\label{eq:compression}
 K_Hu=\sum_{j = 1}^mT_{g_j}H(J_uf_j),\quad u\in Y^*.
\end{equation*}
If \(H\) has finite rank, then \(K_H\) has finite rank and
\begin{equation}\label{eq:compression-bound}
 \norm{K_H}\leq c_{X^*}(\ell_{X^*})\norm{H}.
\end{equation}
In particular, for the identity operator, \(K_{\Id_D}=t_{X^*}(\ell_{X^*})\Id_{Y^*}\).
\end{lemma}
\begin{proof}
For fixed $f_j$, the map $u\mapsto J_uf_j$ is linear and has norm at most $\norm{f_j}$. Thus the definition of $K_H$ gives a bounded linear $K_H:Y^*\to Y^*$ for every $H\in\B(D)$, with the preliminary bound
\[
 \norm{K_H}\leq\norm{H}\sum_{j = 1}^m\norm{g_j}\norm{f_j}.
\]
Assume that $H$ has finite rank. Then
\[
 \ran K_H\subset\sum_{j=1}^m T_{g_j}(\ran H),\quad \text{and} \quad
 \operatorname{rank}K_H\leq m\operatorname{rank}H.
\]
For $u\in Y^*$ and $y\in Y$, the compression
\[
 A_{u,y}=S_yHJ_u:X^* \to X^*
\]
has finite rank and satisfies $\norm{A_{u,y}}\leq\norm{y}\norm{H}\norm{u}$. Therefore,
\begin{align*}
 |(K_Hu)(y)|
 &=\left|\sum_jg_j(S_yH(J_uf_j))\right|\\
 &=|\ell_{X^*}(A_{u,y})|
 \leq c_{X^*}(\ell_{X^*})\norm{H}\norm{u}\norm{y}.
\end{align*}
This immediately yields \eqref{eq:compression-bound}. For $H=\Id_D$, direct computation yields
\[
 S_yJ_uf_j=u(y)f_j,\quad
 (K_{\Id_D}u)(y)=\sum_jg_j(u(y)f_j)=t_{X^*}(\ell_{X^*})u(y).
\]
This implies $K_{\Id_D}=t_{X^*}(\ell_{X^*})\Id_{Y^*}$.
\end{proof}

We proceed to prove Theorem~\ref{thm:main-tensor}.

\begin{proof}[Proof of Theorem~\ref{thm:main-tensor}]
Choose finite tests $$ \ell_{X^*}(A)=\sum_jg_j(Af_j), \quad \ell_{Y^*}(B)=\sum_kb_k(Bu_k), $$
where  $f_j\in X^*$ and $g_j\in X^{**}$, $u_k\in Y^*$ and $b_k\in Y^{**}$, and $$c_{X^*}(\ell_{X^*}),\ c_{Y^*}(\ell_{Y^*})>0, \quad t_{X^*}(\ell_{X^*}) t_{Y^*}(\ell_{Y^*})\neq0. $$
Define $\ell(H)=\ell_{Y^*}(K_H)$ for $H\in\B(D)$. This is a finite test on $\B(D)$ because
\begin{equation}\label{eq:tensor-finite-test}
 \ell(H)=\sum_{j,k}(b_k\circ T_{g_j})\bigl(H(J_{u_k}f_j)\bigr),
 \quad b_k\circ T_{g_j}\in D^*.
\end{equation}
For finite rank $H$, by Lemma~\ref{lem:compression} and the definition of $c_{Y^*}(\ell_{Y^*})$,
\[
 |\ell(H)| \leq c_{Y^*}(\ell_{Y^*})\norm{K_H}\leq c_{Y^*}(\ell_{Y^*})c_{X^*}(\ell_{X^*})\norm{H}.
\]
For the identity operator $\Id_D$, the direct computation in Lemma~\ref{lem:compression} yields
\[
 t_D(\ell) = \ell(\Id_D)=t_{Y^*}(\ell_{Y^*})t_{X^*}(\ell_{X^*}).
\]
Hence $t_D(\ell) \neq 0$. By Lemma~\ref{lem:finite-tests}, $c_D(\ell) \neq 0$. Consequently
$$ 0 < c_D(\ell)\leq c_{Y^*}(\ell_{Y^*})c_{X^*}(\ell_{X^*}). $$
By \eqref{eq:test-beta},
\[
 \beta(D)\geq\frac{|t_D(\ell)|}{c_D(\ell)}
 \geq\frac{|t_{X^*}(\ell_{X^*})|}{c_{X^*}(\ell_{X^*})}\frac{|t_{Y^*}(\ell_{Y^*})|}{c_{Y^*}(\ell_{Y^*})}.
\]
Therefore, \eqref{eq:intro-tensor} follows from Lemma~\ref{lem:finite-tests} by taking all finite tests on $\B(X^*)$ and $\B(Y^*)$.  This finishes the proof.
\end{proof}

\subsection{Tensor amplification}
We first state a well-known proposition on the BAP constant between isomorphic BAP Banach spaces.

\begin{proposition}\label{prop:4.2}
Suppose that $X$ has $\lambda$-BAP and $W : X \to Y$ is an isomorphism. Let $d= \|W\|\|W^{-1}\|$ denote the distortion of $W$. Then $Y$ has $\lambda d$-BAP.
\end{proposition}

We will construct a sequence of dual Banach spaces with exponentially increasing optimal BAP constant via projective tensor products by Theorem~\ref{thm:main-tensor} and Theorem \ref{thm:main-metric}.
% An isomorphism with $\ell_1$ supplies a finite upper bound at each tensor power, which is the additional ingredient required for the direct-sum argument.

\begin{proposition}\label{prop:amplification}
Let $X$ be a Banach space  isomorphic to $\ell_1$ with MAP. Suppose that $ \beta(X^*) > 1$. Assume that $W:X\to\ell_1$ is an isomorphism and $d=\norm{W}\norm{W^{-1}}<\infty$. Define
\[
 X_n=X^{\widehat\otimes_\pi n},\quad E_n=X_n^*,\quad n\in \mathbb{N}.
\]
Then each $X_n$ is separable, isomorphic to $\ell_1$, and has MAP; each $E_n$ is isomorphic to $\ell_\infty$ and has BAP. Moreover,
\begin{equation}\label{eq:general-power-bounds}
 \beta(X^*)^n\leq\beta(E_n)\leq d^n<\infty.
\end{equation}
In particular, $\beta(E_n)$ tend to infinity as $n \to \infty$.
\end{proposition}
\begin{proof}
By Theorem~\ref{thm:main-metric}, such MAP Banach space $X$ does exist. By induction and Theorem~\ref{thm:main-tensor},
\[
\beta(E_n) \geq \beta(X^*)^n.
\]
By iterating the canonical isometric identification
$\ell_1 \widehat\otimes_\pi Y \simeq \ell_1(Y)$ (see
\cite[Example 2.6]{Ryan}), we obtain
\[
\ell_1^{\widehat\otimes_\pi n} \simeq \ell_1(\N^n) \quad \text{isometrically. }
\]

Next, we show that $W^{\otimes n} : X_n \to \ell_1^{\widehat\otimes_\pi n}$ is an isomorphism with distortion $d^n$. By \cite[Proposition 2.3]{Ryan}, both $W^{\otimes n}$ and $(W^{-1})^{\otimes n}$ extend to the projective tensor product spaces, with norms at most $\norm{W}^n$ and $\norm{W^{-1}}^n$, respectively. Their compositions are the identity on the dense algebraic tensor products, and hence on the whole projective tensor products. Therefore, $W^{\otimes n}$ is an isomorphism with distortion at most
$$ \| W^{\otimes n} \| \| (W^{-1})^{\otimes n} \| \leq \| W \|^n \| W^{-1} \|^n = d^n.$$
Since $\N^n$ is countably infinite, $X_n$ is isomorphic to $\ell_1$ and separable. Passing to the adjoint, we obtain an
isomorphism
\[
(W^{\otimes n})^*:
\bigl(\ell_1^{\widehat{\otimes}_\pi n}\bigr)^*
\longrightarrow X_n^*=E_n,
\]
whose distortion satisfies
\[
\begin{aligned}
	\|(W^{\otimes n})^*\|\|((W^{\otimes n})^*)^{-1}\|
	&=\|W^{\otimes n}\|\|(W^{\otimes n})^{-1}\|\\
	&\leq\|W\|^n\|W^{-1}\|^n\\
	&=d^n.
\end{aligned}
\]
Since $\bigl(\ell_1^{\widehat{\otimes}_\pi n}\bigr)^*\simeq \ell_\infty(\mathbb N^n)$ and $\ell_\infty(\mathbb N^n)$ has MAP, Proposition~\ref{prop:4.2} implies
\[
\beta(E_n) \leq d^n.
\]

\iffalse
For a finitely supported tensor $z=\sum_{\mathbf i}c_{\mathbf i}e_{i_1}\otimes\cdots\otimes e_{i_n}$, the definition of the projective norm gives $\norm{z}_\pi\leq\sum_{\mathbf i}|c_{\mathbf i}|$. Choose scalars $\theta_{\mathbf i}$ of modulus at most $1$ such that $\theta_{\mathbf i}c_{\mathbf i}=|c_{\mathbf i}|$. The formula
\[
 B(x^{(1)},\ldots,x^{(n)})
 =\sum_{\mathbf i}\theta_{\mathbf i}x^{(1)}_{i_1}\cdots x^{(n)}_{i_n}
\]
defines an $n$-linear form of norm at most $1$, since the absolute sum is bounded by $\prod_{r=1}^n\norm{x^{(r)}}_1$. Evaluating at $z$ gives the reverse inequality. Finite coordinate tensors are dense in the projective product, and finitely supported arrays are dense in $\ell_1(\N^n)$, so the identification extends onto the completions. Since $\N^n$ is countably infinite, $P_n$ is isomorphic to $\ell_1$, with distortion at most $d^n$, and $D_n$ is isomorphic to $\ell_\infty$ with the same bound. Transport of the usual MAP of $\ell_\infty$ gives $\beta(D_n)\leq d^n$. In particular $\alpha$ is finite.

Associativity of the projective tensor product and Theorem~\ref{thm:main-tensor}, applied successively to $P_{n-1}$ and $X$, give $\beta(D_n)\geq\alpha^n$.
\fi
Finally, since $X$ has MAP, there is a net $(R_\nu)$ of finite rank contractions such that
\[
\| R_\nu x - x \| \to 0, \quad \forall\, x \in X.
\]
% One may direct it by compact sets and positive errors, using inclusion and decreasing errors; hence it approximates every vector.
However, the operators $R_\nu^{\otimes n}$ are still finite rank contractions on $X_n$. For an elementary tensor $x_1\otimes\cdots\otimes x_n$, by the triangle inequality,
\[
 \left\|R_\nu^{\otimes n}(x_1\otimes\cdots\otimes x_n)
       -x_1\otimes\cdots\otimes x_n\right\|_\pi
 \leq\sum_{i=1}^n\norm{R_\nu x_i-x_i}
          \prod_{j\ne i}\norm{x_j} \to 0.
\]
By linearity and density, $R_\nu^{\otimes n}$ converges to the identity operator on $X_n$ pointwisely. This yields that $X_n$ has MAP.
\end{proof}
\subsection{Proof of Theorem~\ref{thm:main-ap}}
We finish the construction by taking an \(\ell_1\)-sum of the tensor-power duals, following the classical argument \cite{FigielJohnson}. Let $G$ be a finitely generated group without Property A. Define
\[
X_n=X_{G,a}^{\widehat\otimes_\pi n},\quad E_n = X_n^*.
\]
By Proposition~\ref{prop:amplification} and Theorem~\ref{thm:3.2}, each $X_n$ is separable, isomorphic to $\ell_1$, and has MAP. Moreover,
\begin{equation}\label{eq:intro-powers}
 (1+b_a)^n\leq\beta(E_n)\leq(1+a)^n.
\end{equation}
In particular, every $E_n$ is isomorphic to $\ell_\infty$ and has BAP, but their optimal BAP constants have no common finite bound. Define
\begin{equation}\label{eq:intro-sum}
 X=\Bigl(\bigoplus_{n\geq1}X_n\Bigr)_{c_0},\quad
 X^*=\Bigl(\bigoplus_{n\geq1}E_n\Bigr)_{\ell_1}
\end{equation}
isometrically under the coordinate pairing. We will show that $X^*$ has AP but fails BAP. We start with the proof of Theorem~\ref{thm:main-ap}.

\

\iffalse
Compact sets in the latter sum have uniformly small tails, so finite-prefix approximation gives AP. Norm-one coordinate compressions would transmit any BAP constant to all $D_n$, contradicting \eqref{eq:intro-powers}. The direct-sum mechanism is classical \cite{FigielJohnson}; the point here is to obtain the unbounded constants among duals with specified MAP preduals. Proposition~\ref{prop:amplification} formulates this tensor-power mechanism independently of the particular group construction.

\begin{lemma}\label{lem:direct-sums}
Let \((P_n)_{n\geq1}\) be nonzero Banach spaces, \(D_n=P_n^*\), and
\[
 V=\Bigl(\bigoplus_{n\geq1}P_n\Bigr)_{c_0}.
\]
Then, under the coordinate pairing,
\begin{equation}\label{eq:canonical-sum}
 V^*=D=\Bigl(\bigoplus_{n\geq1}D_n\Bigr)_{\ell_1}
\end{equation}
isometrically, and
\begin{equation}\label{eq:sum-beta}
 \beta(D)=\sup_{n\geq1}\beta(D_n).
\end{equation}
If every \(D_n\) has AP, then \(D\) has AP. If every \(P_n\) has MAP, then \(V\) has MAP. If every \(P_n\) is separable, then so is \(V\).
\end{lemma}
\fi

\begin{proof}[Proof of Theorem~\ref{thm:main-ap}]
\iffalse
\emph{The isometric dual identification.}
For $(f_n)\in(\bigoplus D_n)_{\ell_1}$ and $(x_n)\in V$, the pairing $\sum_nf_n(x_n)$ is absolutely convergent and has absolute value at most $\sum_n\norm{f_n}\norm{(x_n)}_V$. Conversely, let $f\in V^*$ and let $f_n$ be its restriction to coordinate $P_n$. Fix $N$ and $\eta>0$. For $n\leq N$ choose $x_n\in B_{P_n}$ such that $|f_n(x_n)|>\norm{f_n}-\eta/N$, and change its sign or complex phase to make $f_n(x_n)$ real and nonnegative. The vector supported on these $N$ coordinates has $V$-norm at most $1$, so
\[
 \sum_{n\leq N}\norm{f_n}<\norm{f}+\eta.
\]
Let $\eta\downarrow0$ and then $N\to\infty$. This proves $(f_n)\in(\bigoplus D_n)_{\ell_1}$ with norm at most $\norm{f}$. On finitely supported vectors the two functionals agree, and those vectors are dense in $V$. Thus the representation holds on all of $V$, with equality of norms.
\fi
We divide this proof into three steps.

\

\noindent{\bf Step 1: MAP and separability of $X$.}
Since each $X_n$ is separable, $X$ is also separable. Let $K\subset X$ be any compact subset and $\varepsilon>0$. Since $X$ is the $c_0$-direct sum of $X_n$, for any $x = (x_1, \cdots, x_n, \cdots) \in X$, there is $N > 0$ such that
\[
\sup_{n > N} \| x_n \| < \frac12 \varepsilon.
\]
Take a finite $\varepsilon/2$-net $(x^{(j)})$ of $K$ and choose the common $N$ for this finite net. By the triangle inequality, we obtain
\[
\sup_{x \in K} \sup_{n > N} \| x_n \| < \varepsilon.
\]
Note that the subset $K_n : =\{ x_n : x \in K \}$ is compact in $X_n$. Since each $X_n$ has MAP, there is a sequence of finite rank contractions $S_n \in \F(X_n)$ such that for each $n\in \mathbb{N}$,
\[
\sup_{x_n \in K_n}\| S_n x_n - x_n \| < \varepsilon.
\]
Define $S : X \to X$ by
\[
Sx = (S_1x_1, S_2x_2, \cdots, S_Nx_N, 0, \cdots), \quad \forall\, x \in X.
\]
Then $S$ is a finite rank contraction. Moreover,
\[
\sup_{x \in K} \| Sx - x \| = \sup_{x \in K}\max \left\{ \sup_{n \leq N} \| S_n x_n - x_n \|,  \sup_{n > N} \| x_n \| \right\} < \varepsilon.
\]
This yields that $X$ has MAP.

\

\noindent{\bf Step 2: AP of $X^*$.}
Given compact $H \subset X^*$ and $\varepsilon>0$, by the same argument as in Step 1, there is $N > 0$ such that
\begin{equation}\label{eq:sum-tail}
 \sup_{f \in H}\sum_{n>N}\norm{f_n} < \frac{\varepsilon}{2}.
\end{equation}
The subset $H_n : =\{ x_n : x \in H \}$ is compact in $E_n$.
Since each $E_n$ has BAP, there are $T_n\in\F(E_n)$ for all $n\in\mathbb{N}$ such that
\[
\sup_{f_n \in H_n}\| T_nf_n - f_n \| < \frac{\varepsilon}{2N}
\]
Define $T : X^* \to X^*$ by
\[
 Tf=(T_1f_1,\ldots,T_Nf_N,0,\ldots)
\]
Then $T$ has finite rank and
\[
\| T \| \leq \max_{n \leq N} \| T_n \| < \infty.
\]
This implies $T \in \F(X^*)$. Moreover,
\[
\sup_{f \in H}\| Tf - f \| \leq \sum_{n \leq N} \sup_{f_n \in H_n} \| T_nf_n - f_n \| + \sup_{f \in H}\sum_{n > N} \| f_n \| < \varepsilon.
\]
This yields that $X^*$  has AP.

\

\noindent{\bf Step 3: Failure of BAP of $X^*$.}
We proceed by contradiction. Assume that $X^*$ has $\lambda$-BAP for some $\lambda<\infty$. The coordinate embeddings $i_n: E_n\to X^*$ and projections $p_n: X^* \to E_n$ have norm $1$. For any compact $H_n \subset E_n$ and $\varepsilon > 0$, $i_n(H_n)$ is compact in $X^*$, and hence there is $T \in \F(X^*)$ with $\|T\| \leq \lambda$ such that
\[
\sup_{f \in i_n(H_n)}\| Tf - f \| < \varepsilon.
\]
Therefore, $p_nTi_n \in \F(E_n)$, $\| p_n T i_n \| \leq \lambda$ and
\[
\sup_{f_n \in H_n} \| p_nTi_n (f_n) - f_n \| \leq \sup_{f \in i_n(H_n)}\| Tf - f \| < \varepsilon.
\]
This implies that each $E_n$ has $\lambda$-BAP. Then $\beta(E_n) \leq \lambda$. Consequently,
\[
\sup_{n} \beta(E_n)\leq \lambda<\infty.
\]
However, by Proposition \ref{prop:amplification},
$$  \sup_{n} \beta(E_n)=\infty. $$
This yields a contradiction. Hence, $X^*$ does not have BAP.
\end{proof}

\subsection{Weak MAP without BAP}
We recall the definitions given by \cite[Section 3]{Oja2010}. A Banach space $E$ has weak $\lambda$-BAP, if for every Banach space $F$ and every weakly compact operator $A: E \to F$, there is a net $(T_\alpha)\subset\F(E)$ converging to $\Id_E$ uniformly on norm-compact sets such that
\[
 \limsup_\alpha\norm{AT_\alpha}\leq\lambda\norm{A}.
\]
We say that $E$ has weak BAP if it has weak $\lambda$-BAP for some finite $\lambda\geq1$. Weak MAP means weak $1$-BAP.
% It differs from merely asking for weak or weak-star convergence of the approximating operators: the convergence to the identity here is uniform on norm-compact sets, and the bound concerns composition with each weakly compact $A$.

Lima and Oja showed that AP and weak MAP are equivalent for dual spaces \cite[Corollary 3.4]{LimaOja}; see also \cite[Theorem 3.4]{Oja2010}. Therefore, Theorem~\ref{thm:main-ap} also separates weak MAP from BAP.
\begin{corollary}\label{cor:weak-map}
Over either scalar field, there is a separable Banach space $X$ with MAP such that $X^*$ has weak MAP but does not have BAP.
\end{corollary}

Thus weak BAP does not imply BAP, even in the class of duals with separable MAP preduals. In particular, weak $\lambda$-BAP need not imply $\lambda$-BAP: our dual has weak $1$-BAP but fails $\lambda$-BAP for every finite $\lambda$. This gives a negative answer stronger than a difference between two prescribed finite constants.

\begin{remark}
We point out that our example lies outside the separable-dual \cite{Grothendieck, LT} and RNP \cite{Oja2010} settings in which AP is known to imply MAP. The dual Banach space $E_{G, a}$ is not separable as it is isomorphic to $\ell_\infty$. It embeds as the first coordinate of $X^*$. Hence, $X^*$ is also not separable.

On the other hand, the space $c_0$ fails RNP, and RNP passes to closed subspaces and is preserved by isomorphisms \cite{DU}. Since $E_{G,a}$ contains a closed subspace isomorphic to $c_0$, it fails RNP. Moreover, $\ell_\infty$ contains an isometric copy of $C[0,1]$, and $C[0,1]$ fails Asplund property \cite{NamiokaPhelps1975}. Since Asplund property passes to closed subspaces and is invariant under isomorphisms, $E_{G,a}$ fails Asplund property. Since Asplund property of a Banach space and RNP of its dual are equivalent \cite{NamiokaPhelps1975,Stegall1978}, $E_{G,a}^*$ fails RNP. By coordinate embedding, neither $X^*$ nor $X^{**}$ has RNP.
\end{remark}
\iffalse

\subsection{Quantitative conclusions and limitations}
Theorem~\ref{thm:bap-bounds} supplies an explicit obstruction on $|S|+3$ label indicators, expressed through $\kappa_S(G)$. It does not compute $\kappa_S(G)$ for Osajda's group or the exact value of $\beta(E_{G,a})$. Neither endpoint of the estimate is asserted to be sharp. Corollary~\ref{thm:group} proves only
\[
 E_{G,a}\text{ has MAP}\quad\Longrightarrow\quad G\text{ has Property A};
\]
the converse is not established here. Similarly, Theorem~\ref{thm:main-tensor} is a lower bound; it does not assert equality or a finite upper bound for arbitrary tensor-product duals. In Proposition~\ref{prop:amplification}, the upper bound comes specifically from the isomorphism with $\ell_1$.

There is a general reason that an unbounded family of dual approximation constants is enough for \ref{q:ap-bap}. If every dual with AP had BAP but no common finite bound on its optimal constants existed, choose duals $P_n^*$ with AP and $\beta(P_n^*)>n$. Lemma~\ref{lem:direct-sums} would give a dual $\ell_1$-sum with AP and no BAP, a contradiction. Our construction realizes this mechanism explicitly: a strict constant gap is multiplied by tensor powers, and the $c_0$-sum preserves separability and MAP on the predual. This is the link between the isometric separation in \ref{q:bap-map} and the isomorphic separation in \ref{q:ap-bap}.
\fi
\bigskip

\noindent{\bf AI Statement.}
During the preparation of this paper, the authors used AI tools to assist in developing the constructions and improve the language and exposition. The idea of using Property A in the proof emerged through discussions with GPT-6 Astra. Suggestions from these AI interactions were subsequently examined, reformulated, incorporated into the manuscript and independently verified by the authors. The authors take full responsibility for all mathematical content in the final manuscript.

%\noindent{\bf Acknowledgments.}

\enlargethispage{3\baselineskip}

\end{document}